\documentclass{degm}
\usepackage{amsmath,amssymb,amstext}

\evensidemargin\oddsidemargin
\newcommand{\skipit}[1]{\relax}

\newtheorem{Theo}{Theorem}
\newcommand{\Tb}{\begin{Theo}}
\newcommand{\Te}{\end{Theo}}
\newcommand{\proof}{\it Proof. ~\rm}

\newtheorem{Cor}{Corollary}[Theo]
\newcommand{\Cb}{\begin{Cor}}
\newcommand{\Ce}{\end{Cor}}
\newtheorem{Lem}{Lemma}
\newcommand{\Lb}{\begin{Lem}}
\newcommand{\Le}{\end{Lem}}
\def\D{\displaystyle}
\def\T{\textstyle}
\def\fn{\footnotesize}

\newcommand{\Coro}[2] {\bigskip\noindent{\bf Corollary #1~~}  {\it
#2}\bigskip}

\def\klfrac#1#2{\mathsurround0pt\T{%
   \frac{\hspace{0.03em}\D%
   \rule[-.7ex]{0pt}{1.5ex}%
   \mbox{\fn$#1$}\hspace{0.03em}}%
   {\hspace{0.03em}\D \mbox{\fn$#2$}%
   \rule[.2ex]{0pt}{1ex}%
   \vphantom{\overline{\mbox{\fn$#2$}}}%
   \hspace{0.03em}}}}

\newcommand{\nn}{\nonumber}

\newcommand{\ga}{\gamma}
\newcommand{\sig}{\sigma}
\newcommand{\de}{\delta}
\newcommand{\ka}{\kappa}

\def\ve{\varepsilon}

\newcommand{\ld}{\ldots}

\newcommand{\ca}[1]{{\cal #1}}

\newcommand{\setN}{\mathbb N}

\newcommand{\meq}{\;=\;}
\newcommand{\defeq}{\;:=\;}
\newcommand{\mge}{\;\ge\;}
\newcommand{\mg}{\;>\;}
\newcommand{\ml}{\;<\;}
\newcommand{\mle}{\;\le\;}

\newcommand{\myin}{\!\in\!}
\def\impl{\Rightarrow} 

\newcommand{\astA}{\mathop{*_{\!\scriptscriptstyle A}}}

\newcommand{\udiv}{\mathop{|\hspace{-1pt}|}}

\newfont{\syten}  {cmsy10  scaled 1000}
\newfont{\syeight}{cmsy8   scaled 1000}
\newfont{\sysix}  {cmsy6   scaled 1000}
\newcommand{\BOX}{\mbox{\syten \char"74\hspace{-.6625em}\char"75}}
\newcommand{\Dqed}{\tag*{\BOX}}
\renewcommand{\qed}{\hfill\BOX}

\newcommand{\comment}[1]{{\it\scriptsize{#1}\vspace{-.5\baselineskip}

}}
\def\ttg{\spa|\spa}
\def\ttk{\mathrel{|}}
\def\tt{\mathchoice{\ttg}{\ttg}{\ttk}{\ttk}}
\newcommand{\spa}{\hspace{1pt}}

\begin{document}
\renewcommand{\today}{April, 2003}

\title{\vspace*{10mm}The maximal order of a class of multiplicative arithmetical
 functions}
\author{L\'aszl\'o T\'oth \thanks{Supported partially by the
Hungarian National Foundation for Scientific Research under grant
OTKA T031877} (P\'ecs) and Eduard Wirsing (Ulm)\\[1mm]\\
Annales Univ. Sci. Budapest., Sect. Comp., \\ {\bf 22} (2003), 353-364}
\maketitle

\markboth{\hfill L\'aszl\'o T\'oth and Eduard Wirsing \hfill}
{\hfill Maximal order of a class of multiplicative
 functions \hfill}
\vspace{-3\baselineskip}
\thispagestyle{empty}

\begin{abstract} We prove simple theorems concerning the maximal
order of a large class of multiplicative functions. As an
application, we determine the maximal orders of certain functions of
the type $\sigma\!_A(n)= \sum_{d\in A(n)} d$, where $A(n)$ is a
subset of the set of all positive divisors of $n$, including the
divisor-sum function
$\sigma(n)$ and its unitary and exponential analogues.
We also give the minimal order of a new class of Euler-type
functions, including the Euler-function $\phi(n)$ and its unitary
analogue.
\end{abstract}



\section{Introduction}
Let $\sigma(n)$ and $\phi(n)$ denote, as usual, the sum of all positive
divisors of $n$ and the Euler function, respectively. It is well-known, that
\begin{gather}
\limsup_{n\to \infty} \frac{\sigma(n)}{n\log \log n}\meq
e^{\gamma}\,,\label{Sig} \\[1ex]
\liminf_{n\to \infty} \frac{\phi(n)\log \log n}{n}
\meq e^{-\gamma}\,,\label{Phi}
\end{gather}
where $\gamma$ is Euler's constant. These results go back to the work of T.H.\
Gronwall \cite{Gr13} and E.\ Landau \cite{La09} and have been established for a
number of modified $\sig$\!- and $\phi$\!- functions.

One such modification relates to unitary divisors $d$ of $n$, notation
$d\udiv n$, meaning that $d\tt n$ and $(d,n/d) = 1$. The corresponding
$\sig$\!- and $\phi$\!-functions are defined by $\sig^*(n)=\sum_{d\udiv n} d$
and $\phi^*(n)=\#\{1\le k\le n\,; (k,n)_*=1\}$, where $(k,n)_*$ denotes the
largest divisor of $k$ which is a unitary divisor of $n$. These functions are
multiplicative and for prime powers $p^\nu$ given by
$\sig^*(p^\nu) =p^\nu+1,\; \phi^*(p^\nu)= p^\nu-1$, see \cite{Co60,McC86}.
They are treated, along with other multiplicative functions, in \cite{ChSi85}
with the result that
\begin{align}\label{unitary Sig}
\limsup_{n\to \infty} \frac{\sigma^*(n)}{n\log\log n}
&\meq\frac{6}{\pi^2}e^{\gamma}\,,
\end{align}
while $\phi^*$ gives again \eqref{Phi}.
(Actually \eqref{unitary Sig} is  written incorrectly in \cite{ChSi85}
with the factor $6/\pi^2$ missing).

In \cite{FaSu89} it is shown that \eqref{unitary Sig} holds also for
$\sigma^{(e)}(n)$\!, the sum of exponential divisors of $n.$
(A number
$d=\prod p^{\de_p}$ is called an exponential divisor of $n=\prod p^{\nu_p}$
if $\de_p|\nu_p$ for all $p.)$

These and a number of similar results from literature refer to rather
special functions. Textbooks dealing with the extremal order of
arithmetic functions also treat only particular cases, see
\cite{HaWr60,Ap76,Te95}.
It should be mentioned that a useful result concerning the maximal order of a
class of prime-independent functions, including the number of all divisors,
unitary divisors and exponential divisors, is proved in \cite{SuSi75}.
\skipit{\comment{They prove the following: Let $f$ be
positive and satisfying
$f(n)=O(n^\beta)$ for some fixed $\beta >0$. Let $F$ be
multiplicative with $F(p^\nu)=f(\nu)$. Then
$\limsup_{n\to \infty} \frac{\log F(n) \log \log n}{\log n} = \sup_{m}
\frac{\log f(m)}{m}$.}}

In the present paper we develop easily applicable theorems for determining
\begin{equation*}
L\meq L(f)\defeq \limsup_{n\to\infty}\frac{f(n)}{\log\log n}
\end{equation*}
where $f$ are nonnegative real-valued multiplicative functions.
Essential parameters are
\begin{equation*}
\rho(p)\meq
\rho(f,p)\defeq\sup_{\nu\ge0}f(p^\nu)\,,
\end{equation*}
for the primes $p,$ and the product
\begin{equation*}
R\meq R(f)\defeq\prod_p \Big(1-\frac1p\Big)\rho(p)\,.
\end{equation*}
These theorems can, in particular, be used to obtain the maximal or minimal
order, respectively, of  generalized $\sigma$\!- and $\phi$\!-functions which
arise in connection with Narkiewicz-convolutions of arithmetic functions.

\section{General results}
We formulate the conditions for lower and upper estimates for $L$ separately.
Note that $\rho(p)\ge f(p^0)=1$ for all $p.$

\Tb\label{upper}
Suppose that $\rho(p)<\infty$ for all primes $p$ and that the
product $R$ converges unconditionally (i.e.\ irrespectively of order), improper
limits being allowed, then
\begin{equation}\label{res}
L\mle e^\gamma R\,.
\end{equation}
\Te

A different assumption uses\vspace{-1ex}
\Tb\label{upper2}
Suppose that $\rho(p)<\infty$ for all $p$ and that the product $R$ converges,
improper limits being allowed, and that
\begin{equation}\label{altern}
\rho(p)\mle 1 + o\Big(\frac{\log p}{p}\Big)\,.
\end{equation}
then \eqref{res} holds.
\Te

{\bf Remark.} Neither does \eqref{altern} plus convergence imply unconditional
convergence nor vice versa.
\bigskip

To establish $e^\ga R$ also as the lower limit more information is
required: The suprema $\rho(p)$ must be sufficiently well approximated
by not too large powers \mbox{of $p.\!$}
\Tb\label{lower}
Suppose that $\rho(p)<\infty$ for all primes $p,$ that for each prime
$p$
there is an exponent $e_p=p^{o(1)}\in\setN$ such that
\begin{equation}\label{suitab}
\prod_p f(p^{e_p})\rho(p)^{-1}\mg 0\,,
\end{equation}
and that the product $R$ converges, improper limits being allowed.
Then
\begin{equation*}
L\mge e^\gamma R\,.
\end{equation*}
\Te

\Coro{1.}{If for all $p$ we have $\rho(p)\le(1-1/p)^{-1}$
and there are $e_p$
such that $f(p^{e_p})\ge 1+1/p$ then
\begin{equation*}
L\meq e^\gamma R\,.
\end{equation*}
In other words: The maximal order of $f(n)$ is $e^\ga R\log\log n$.}
\medskip

Formally $R$ becomes infinite if there is a nonempty set $\ca S$ of primes for
which $\rho(p)=\infty$. So one might expect that the assumptions of
Theorem\,\ref{lower} taken for all $p$ with finite $\rho(p)$ would imply
$L=\infty$. Surprisingly enough this is true only for rather thin sets $\ca S.$
But note that for $p\in \ca S$ there is no substitute for the $f(p^{e_p})$
approximating $\rho(p)$.

We begin by stating what the above theorems imply if one ignores the numbers
with prime factors from a given set $\ca S$ of primes. For any such set define
\begin{gather*}
N(\ca S)\defeq\{n\!: n\in\setN,\;p\tt n\impl p\in\ca S\}\,,\quad
C(\ca S)\defeq \{n\!: n\in\setN,\;p\tt n\impl p\notin \ca S\}\,.
\end{gather*}

\Coro{2.}{Modify the assumptions of Theorems 1, 2, and 3 by
replacing $R$ with
\begin{equation*}
R_{\ca S}\meq R_{\ca S}(f) \defeq \prod_{p\notin \ca S}
\Big(1-\frac1p\Big)\rho(p)\,,
\end{equation*}
$L$ with
\begin{equation*}
L_{\ca S}\meq L_{\ca S}(f)\defeq\limsup_{n\to\infty,\,n\in C(\ca S)}
\,\frac{f(n)}{\log\log n}\,,
\end{equation*}
condition \eqref{altern} with
\begin{equation}\label{alternS}
\rho(p)\mle 1 + o\Big(\frac{\log p}{p}\Big)\quad\text{for}
\quad p\notin\ca S\,.
\end{equation}
and \eqref{suitab} with
\begin{equation}\label{suitabS}
\prod_{p\notin \ca S}f(p^{e_p})\rho(p)^{-1}\mg 0\,.
\end{equation}
Assume further that
\begin{equation*}
\sum_{p\in\ca S}\frac1p\ml \infty\,.
\end{equation*}
Then
\begin{equation*}
L_{\ca S}\mle e^\ga\prod_{p\in\ca S}\Big(1-\frac1p\Big)\cdot
 R_{\ca S}\,,\quad
L_{\ca S}\mge e^\ga\prod_{p\in\ca S}\Big(1-\frac1p\Big)\cdot
 R_{\ca S}\,,
\end{equation*}
respectively. This applies even if $\rho(p)=\infty$ for some or all
of the $p\in \ca S.$}

\Tb\label{infin}
Let $\ca S$ be a set of primes such that
\begin{equation}\label{thin}
\sum_{p\in\ca S}\frac1p\ml \infty\,.
\end{equation}
If $\rho(p)=\infty$ exactly for the $p\in \ca S$, if \eqref{suitabS}
holds and $R_{\ca S}>0$, then $L=\infty$.
Condition \eqref{thin} must not be waived.
\Te
In fact there are counter-examples for any set $\ca S$ for which
$\sum 1/p$
diverges.

\section{The proofs}
{\it Proof of Theorem\,\ref{upper}.}
An arbitrary $n=\prod p^{\nu_p}$ we write as $n=n_1n_2$ with $ n_1:=
\prod_{p\le\log n}p ^{\nu_p}$. Mertens's formula $\prod_{p\le x}(1-1/p)^{-1}
\sim e^\ga\log x$ and\vspace{2pt} the definition of $\rho(p)$ imply
\begin{align}
f(n_1)&\meq\prod_{p\le\log n}\!f(p^{\nu_p})\mle\prod_{p\le\log n}\rho(p)\nn\\
&\meq\prod_{p\le\log n}\!\!\Big(1-\frac1p\Big)^{-1}\cdot
\prod_{p\le\log n}\!\!\Big(1-\frac1p\Big)\rho(p)\,,\nn\\[1ex]
f(n_1)&\mle\big(1+o(1)\big)e^\ga R\,\log\log n\,
\quad \text{ as } n\to\infty\,.\label{n1}
\end{align}
Let $a$ denote the number of prime divisors in $n_2$. Then $a\le\log n
/\log\log n$. There is nothing to prove if $R=\infty$, so let $R<\infty$.
Using the unconditional convergence
\begin{align}
f(n_2)&\mle\prod_{p|n,\,p>\log n}\!\Big(1-
\frac1p\Big)\rho(p)\,\cdot\!\!
\prod_{p|n,\,p>\log n}\!\Big(1-\frac1p\Big)^{-1}\nn\\
&\mle\big(1+o(1)\big)\cdot\Big(1-\frac1{\log n}\Big)^{-a}\nn\\
&\meq\big(1+o(1)\big)\,e^{O(1/\log\log n)}\;\to\;1\,.\label{n2}
\end{align}
Combining \eqref{n1} and \eqref{n2} finishes the proof. \qed

{\it Proof of Theorem\,\ref{upper2}.\/} There is no change in the
estimation of $f(n_1)$.
For $n_2$ we have
\begin{align*}
f(n_2)&\mle \Big(1+o\big(\klfrac{\log\log n}{\log
n}\big)\Big)^\frac{\log n}{\log\log n}\\
&\meq 1+ o(1)\,. \Dqed
\end{align*}

{\it Proof of Theorem\,\ref{lower}.\/} We treat the case of
 proper convergence
only. There is nothing to prove if $R=0$ and the changes for
$R=\infty$ are
obvious.  For given $\ve$ take $P$ so large that
\begin{equation}\label{P l 1}
\prod_{p>P}f(p^{e_p})\rho(p)^{-1}\mge 1-\ve
\end{equation}
and choose exponents
$k_p$ for the $p\le P$ such that
\begin{equation}\label{P l 2}
\prod_{p\le P}f(p^{k_p})\mge (1-\ve)\prod_{p\le P}\rho(p)\,.
\end{equation}
Keeping $P$ and the $k_p$ fixed let $x$ tend to infinity and consider
\begin{equation*}
n(x)\;:=\;\prod_{p\le P}p^{k
_p}\prod_{P<p\le x}p^{e_p}\,.
\end{equation*}
Now on the one hand, using \eqref{P l 1} and \eqref{P l 2}
we see
\begin{align*}
f\big(n(x)\big)\prod_{p\le x}\!\Big(1-\frac1p\Big)
&\mge(1-\ve)\prod_{p\le x}\!\Big(1-\frac1p\Big)\rho(p)\cdot\!
\prod_{P<p\le x
}f(p^{e_p})\rho(p)^{-1}\nonumber\\
&\mge(1-\ve)^2\big(1+o(1)\big) R
\end{align*}
and with Mertens's formula again
\begin{equation}\label{fnx}
f\big(n(x)\big)\mge (1-\ve)^2\big(1+o(1)\big)R\,e^\ga\log x\,.
\end{equation}
On the other hand, since $e_p=p^{o(1)},$ we have
$$
\log n(x)\mle\sum_{p\le P}k_p\log p\; +\! \sum_{P<p\le x}e_p\log p\mle
x^{o(1)}\sum_{p\le x}\log p\meq x^{1+o(1)}\,,
$$
and therefore
$$
\log\log n(x)\mle \big(1+o(1)\big)\log x\,.
$$
Together with \eqref{fnx} this yields the lower bound

\begin{equation*}
\limsup_{x\to\infty}\,\frac{f\big(n(x)\big)}{\log\log n(x)}\mge
(1-\ve)^2 R\,e^\gamma
\end{equation*}
with arbitrary $\ve>0$. \qed

{\it Proof of Corollary 1.\/} Apply Theorems\, \ref{upper} (or
\ref{upper2}) and \ref{lower}. \qed

{\it Proof of Corollary 2.\/} To see this one applies the theorems to
the multiplicative function $f^*$
defined by $f^*(n)=f(n)$ for $n\in C(\ca S)$ and $f(n)=1$ for
$n\in N(\ca S)$. One finds
$$
L(f^*)=L_{\ca S}(f)\,,\quad
R(f^*)=R_{\ca S}(f)\prod_{p\in\ca S}\Big(1-\frac1p\Big)\,,
$$
and \eqref{suitabS} implies \eqref{suitab} for
$f^*$ because $\prod_{p\in\ca S}(1-1/p)$ converges absolutely.

Note also that for any sequence of numbers $n=n_1n_2$ tending to $\infty$,
where $n_1\in N(\ca S)\,, \;n_2\in C(\ca S)$, we have $f^*(n)/\log\log n =
f(n_2)/\log\log(n_1n_2)$, hence  $\limsup f^*(n)/\log\log n=0$ if $n_2$ stays
bounded, and $\le\limsup_{n_2}f(n_2)/\log\log n_2$ otherwise, with equality if
$n_1$ is bounded. Thus $L(f^*)=L_{\ca S}$. \qed

{\it Proof of Theorem\,\ref{infin}.\/}
{\bf ~I.} Assume \eqref{thin}.
With any $n_1\in N(\ca S)$ we have
\begin{equation*}
L\mge\limsup_{n_2\in C(\ca S)}\,\frac{f(n_1)f(n_2)}{\log\log (n_1n_2)}
\meq f(n_
1)L_{\ca S}\,.
\end{equation*}
 From Corollary\,2, as it refers to Theorem\,\ref{lower}, we have
$L_{\ca S}>0$
and $f(n_1)$ can be chosen arbitrarily big.

{\bf II.} Assume that \eqref{thin} does not hold. We shall construct a
counter-example.
The assumption implies that
\begin{equation*}
g(x)\defeq\prod_{p\in\ca S,\,p\le x}\Big(1+\frac1p \Big)
\end{equation*}
tends to $\infty$ as $x\to\infty$. Choose an increasing sequence of
numbers
$q_j=p_j^{\nu_j}$ with $p_j\in\ca S$ and $\nu_j$ so large that $g(\log
q_j)\ge
j^j$ for all $j,$ and such that every prime $p\in \ca S$ occurs
infinitely
often in the sequence of the $p_j$. Put $f(q_j)=j$ for all
$j\in\setN$ and
$f(p^\nu)=1+1/p$ for all $p^\nu$ that are not among the $q_j$. Then,
obviously,
$\rho(p)=\infty$
for $p\in\ca S$ and $\rho(p)=1+1/p$ for $p\notin \ca S$. The
product
\begin{equation*}
R_{\ca S}\meq\prod_{p\notin\ca S}\Big(1-
\frac1p\Big)\Big(1+\frac1p\Big)
\end{equation*}
converges absolutely and so does (choosing $e_p=1$)
$\prod_{p\notin\ca S}f(p^1)
/\rho(p)=1$.  Any
$n\in\setN$ can be written as $n=n_1n_2$ where $n_1$ collects from the
canonical representation of $n$ those prime powers that occur among
the $q_j$
while the rest compose $n_2$. For given $n$ let $k:=\max\{j;\,q_j
\udiv n_1\}$.
Then $f(n_1)\le k!=o(k^k)=o\big(g(\log q_k)\big)=o\big(g(\log
n)\big)$ by construction.
Now for any $n\in\setN$
\begin{align*}
f(n)&\meq f(n_1)f(n_2)
\meq \prod_{p\tt n,\,p\notin\ca S}\Big(1+\frac1p\Big)\cdot
o\big(g(\log n)\big)\\
&\meq o\bigg(\prod_{p\le \log n}\Big(1+\frac1p\Big)\bigg)\cdot\D
\prod_{p\tt n,\,p\ge\log n}\Big(1+\frac1p\Big)\\
&\mle o(\log\log n)\cdot
\Big(1+\frac1{\log n}\Big)^\frac{\log n}{\log\log n}
\meq o(\log\log n)\,,
\end{align*}
hence $L=0$.\qed

\section{Applications}

A general frame for generalizations of the $\sigma$\!- and $\phi$\!-functions
mentioned in the introduction can be found in Narkiewicz \cite{Nar63}. Assume
that for each $n$ a set $A(n)$ of divisors of $n$ is given and consider the
$A$\!-convolution $\astA$ defined by
\begin{equation}\label{Nark convol}
\big(f\astA g\big)(n)\;:=\;\sum_{d\in A(n)}\!\!f(d)\,g\!\left(\frac
nd\right).
\end{equation}

Properties of convolution \eqref{Nark convol} and of arithmetical functions
related to it have been studied extensively in the literature, see
\cite{Nar63, McC86}. The system $A$ is called multiplicative if $A(n_1
n_2)=A(n_1)A(n_2)$ for coprime $n_1,n_2,$ with elementwise multiplication of
the sets, and not all $A(n)$ empty. Such a divisor system can be described
by the sets
$A\hspace{-.6pt}E_p(\nu)$ of admissible exponents,
$$
A\hspace{-.6pt}E_p(\nu)\;:=\;\{\delta;\ p^\delta\myin A(p^\nu)\}\,.
$$
The $A$\!-convolution of any two multiplicative functions $f$ and $g$
is
multiplicative if and only if $A$ is multiplicative. In particular
multiplicativity of $A$ implies multiplicativity of the modified
divisor function
$$
\sig\!_A(n)\;:=\;\sum_{d\in A(n)}d\,.
$$

As a natural means to define an Euler-function attached to $A$ we consider the relation
\begin{equation}\label{phiA}
\sum_{d\in A(n)}\phi_A(d) \meq n\, , \quad n\ge 1.
\end{equation}
This need not be solvable; there is, however, the following
\Tb\label{sol phi}
If the divisor system $A$ is multiplicative then \eqref{phiA} has a solution
if and only if $n\in A(n)$ for all $n\in\setN$. In this case  the solution
$\phi\!_A$ is unique and is a multiplicative function with $1\le \phi\!_A(n)
\le n$ for all $n\in\setN$.
\Te
\proof Suppose a solution exists. Then by induction on $\nu$ the recursion
\begin{equation}\label{phi exp}
\sum_{\delta\in A\hspace{-.6pt}E_p(\nu)}\phi_A(p^\delta) \meq p^\nu
\end{equation}
implies  that $1\le\phi_A(p^\nu)\le p^\nu$ and (therefore) $\nu\in
A\hspace{-.6pt}E_p(\nu)$ for all $\nu\!:\;p^\nu\in A(p^\nu)$.  It follows from
the multiplicativity of $A$ that $n\in A(n)$ for all $n.$ If, on the other
hand, $n\in A(n)$ for all $n,$ then \eqref{phi exp} can be solved recursively
and the multiplicative function defined from the $\phi_A(p^\nu)$ solves
\eqref{phiA}. This is in fact the only solution since $\phi\!_A(n)=
n-\sum_{d\in A(n)\smallsetminus\{n\}}\phi\!_A(d)$. \qed
\bigskip

With suitable additional conditions on $A$ we give the maximal and
minimal orders of $\sigma\!_A$ and $\phi\!_A$, respectively.
Extremal orders of such functions have not been investigated in
the literature.

Obviously $\sig\!_A(n)\le \sig(n)$ and if for any $\nu$ we have
$p^\nu,\,p^{\nu-1}\myin A(p^\nu)$ then $\sig \!_A(p^\nu)\ge p^\nu+p^{\nu-1}.$
So Corollary\,1 applies to $f(n)=\sig\!_A(n)/n$ and gives

\Tb\label{sigmaA}
Let the system $A$ of divisors be multiplicative and suppose that for
each prime $p$ there is an exponent $e_p$ such that
$$
p^{e_p},\, p^{e_p-1} \in A(p^{e_p})
$$
and $e_p = p^{o(1)}.$ Then
$$
\limsup_{n\to\infty} \frac{\sigma\!_A(n)}{n \log \log n} \meq e^\gamma
\prod_p \Big(1-\frac1p\Big) \sup_{\nu \ge 0}
\frac{\sigma\!_A(p^{\nu})}{p^{\nu}}\,,
$$
where the product converges.
\Te

{\bf Remarks.} The quotients $\sigma(p^\nu)/p^\nu$ are of the form
$\sum \ve_ip^{-i},\ \ve_i\in\{0,1\}$, and the set of such numbers is compact.
Therefore each $\sup_\nu \big(\sig\!_A(p^\nu) p^{-\nu}\big)$ is itself of
this form and we have for each prime $p$ a finite or infinite
sequence of exponents $a_i$ such that $2\le a_1< a_2<\ldots$ and
$$
\Big(1-\frac1p\Big) \sup_\nu \frac{\sigma\!_A(p^{\nu})}{p^{\nu}}
\meq 1-\frac1{p^{a_1}}+\frac1{p^{a_2}}-\frac1{p^{a_3}}+\ldots\;.
$$

The formulae \eqref{Sig} and \eqref{unitary Sig} are obvious consequences of
Theorem \ref{sigmaA}. In the standard case $e_p$ is arbitrary, we have
$(1-1/p)\rho(p)=1$ for all $p,$ hence $R=1$. With unitary and exponential
divisors the only admissible choices are $e_p=1$ and $e_p=2$, respectively and
$(1-1/p)\rho(p)=1-1/p^2,$, hence $R=\zeta(2)^{-1}=6/\pi^2$ in both cases.
\bigskip

We turn to $\phi\!_A,$ assuming again that $A$ is multiplicative and, in view
of Theorem \ref{sol phi}, that always $\nu\in A\hspace{-.6pt}E_p(\nu)$.
In order to determine the minimal order of $\phi\!_A$ consider the function $f(n):=n/\phi\!_A(n)$.

For all $p$ and $\nu\ge 1$ we have
$\phi\!_A(p^\nu)\ge p^\nu-\phi\!_A(p^{\nu-1})-\ldots-\phi\!_A(1)\ge
p^\nu-p^{\nu-1}-\ldots-1$, which gives
$$
f(p^\nu)\ml \frac{p-1}{p-2}\,,\qquad\rho(p)\mle \frac{p-1}{p-2}\,.
$$
Note that $\rho(2)$ may equal $\infty$.
If moreover $e-1\in A\hspace{-.6pt}E_p(e)$ for some $e=e_p\ge 1$ then, on the
other hand, $\phi_A(p^e)\le p^e -\phi_A(p^{e-1})\le
p^e-p^{e-1}+p^{e-2}+\ldots+1$ if $e\ge 2$, and $\phi_A(p)\le p-1$ if $e=1$.
Therefore
\begin{gather*}
f(p^e)\mge \frac{p(p-1)}{p^2-2p+2}\,,\\[1ex]
f(p^e)\rho(p)^{-1}\mge \frac{p(p-2)}{p^2-2p+2}\meq 1-\frac{2}{p^2-2p+2}\,,
\end{gather*}
which is positive and yields a convergent product for $p\ge 3$.

Note that for  powers of $2$ there is no non-trivial lower estimate for
$\phi_A(n)/n$ whithout further conditions on $A.$ This is shown by the
following example.  Let $\ca N=\{n_1,n_2,\ldots\}\subset\setN$, $n_1<n_2<\ld$,
and put $A\hspace{-.6pt}E_2(n):= \{0,1,\ldots,n\}$ for $n\myin \ca N$ and
$A\hspace{-.6pt}E_2(n):= \{n\}$ for $n\notin\ca N$.  Then the recursion gives
$\phi_A(2^n)=2^n$ for $n\notin\ca N$ but $\phi_A(2^{n_j})= 2^{n_{j-1}}$ for the
$n\myin\ca N$, where $n_0=0$.  Hence it is possible to have
$\rho(2)=\sup_{\nu} 2^\nu/\phi\!_A(2^\nu) = 2^{\sup_j (n_j-n_{j-1})}=\infty$.

Thus applying Corollary 1 or Theorem \ref{infin} with $\ca S=\{2\}$ we
obtain

\Tb \label{phi}Let $A$ be multiplicative and $n\myin A(n)$ for all $n.$
Assume that
for each prime $p>2$ there is an exponent $e_p$ such that $p^{e_p-
1}\myin A(p^{e_p})$ and $e_p = p^{o(1)}.$  Then
$$
\liminf_{n\to\infty}\frac{\phi_A(n) \log \log n}n \meq  e^{-\gamma}
\prod_p\Big(1-\frac1p\Big)^{-1} \inf_\nu
\frac{\phi_A(p^{\nu})}{p^{\nu}}\,.
$$
The product converges for $p>2$; the first factor may vanish.
\Te

For the standard Euler function $\phi(n)$ and for its unitary
analogue $\phi^*(n)$ we regain \eqref{Phi}.

For the system of exponential divisors one has $\phi_A(1)=1$ because
of
multiplicativity. The recursion
$\sum_{\ka|\nu}\phi_A(p^\ka)=p^\nu$ is
solved by $\phi_A(p^\nu)=\sum_{\ka|\nu}\mu(\nu/\ka)p^\ka.$ Again the
minimum of
$\phi_A(p^\nu)/p^\nu$ is $1-1/p$, it is taken for $\nu=e_p=2$ and once
more \eqref{Phi} follows.

\bigskip

\noindent
\parbox[t]{.47\textwidth}{\small L\'aszl\'o T\'oth\\
University of P\'ecs\\
Institute of Mathematics and Informatics\\
Ifj\'us\'ag u. 6,\\
7624 P\'ecs, Hungary\\
ltoth@ttk.pte.hu}
\hfill
\parbox[t]{.38\textwidth}{\small Eduard Wirsing\\
Universit\"at Ulm\\
Helmholtzstra{\ss}e 22,\\
D--89069 Ulm, Germany\\
Tel.: +49--731--50--23565\\
wirsing@mathematik.uni-ulm.de}

\end{document}